\documentclass[coll]{imsart}
\RequirePackage[OT1]{fontenc}
\usepackage{amsthm,amsmath}

\usepackage{amsfonts,latexsym,amssymb,mymathbbol}
\usepackage{graphics,graphicx}

\pubyear{0000}
\volume{0}
\volumetitle{Volume Title}
\firstpage{1}
\lastpage{6}
\arxiv{math.PR/0000000}

\startlocaldefs
\def\:{:\,}

\def\sec{\setcounter{equation}{0}\setcounter{figure}{0}}

\def\EC{Euler characteristic\ }

\def\bb{\mathbb }

\def\cal{\mathcal}

\def\LK{Lipschitz-Killing}

\def\P{{\bb P}}
\def\E{{\bb E}}

\def\Z{{\bb Z}}

\def\e{{\varepsilon}}

\def\real{{\bb{R}}}

\def\Rk{\mathbb{R}^k}

\def\definedas{\stackrel{\Delta}{=}}

\def\Min{{\cal M}}

\newcommand{\sqbinom}[2]{\begin{bmatrix}#1 \\ #2 \end{bmatrix}}

\def\qed{\hfill $\Box$\newline}

\newcommand{\bc}{\begin{center}}
\newcommand{\ec}{\end{center}}

\newcommand{\eex}{\end{exercise}}
\newcommand{\bex}{\begin{exercise}}

\newcommand{\smc}{\scshape}

\newcommand{\sign}{\text{\rm sgn}}

\newcommand{\Tr}{{\rm Tr}}

\newcommand{\lips}{{\cal L}}

\def\text{\mbox}

\newtheorem{Theorem}{Theorem}[section]

\newtheorem{Lemma}[Theorem]{Lemma}

\newtheorem{Corollary}[Theorem]{Corollary}

\newtheorem{definition}[Theorem]{Definition}
\newcommand{\beq}{\begin{eqnarray}}
\newcommand{\eeq}{\end{eqnarray}}
\newcommand{\beqq}{\begin{eqnarray*}}
\newcommand{\eeqq}{\end{eqnarray*}}

\DeclareMathOperator{\rk}{rank}

\DeclareMathOperator{\lf}{lf}

\def\R{\mathbb{R}}
\def\Z{\mathbb{Z}}

\def\e{\varepsilon}

\def\vp{\varphi}

\newcommand{\mean}[1] {\E\left\{{#1}\right\}}

\newcommand{\indic}{{\mymathbb{1}}} 

\newcommand{\LL}{{\mathcal{L}}}
\newcommand{\MM}{{\mathcal{M}}}

\newcommand{\set}[1]{\left\{#1\right\}}
\newcommand{\sbrk}[1]{\left[#1\right]}

\newcommand{\param}[1]{\left(#1\right)}
\newcommand{\abs}[1] {\left| {#1}\right|}

\newcommand{\iprod}[1] {\left< {#1}\right>}
\newcommand{\cP}{{\mathcal{P}}}

\newcommand{\cech}{\v{C}ech }

\newcommand{\EI}{{Euler integral}}

\newcommand{\proofsmodeb}[1]{{aaa}}

\numberwithin{equation}{section}
\theoremstyle{plain}

\endlocaldefs

\begin{document}

\begin{frontmatter}
\title{Persistent Homology for Random Fields and Complexes}
\runtitle{Persistent Homology for Random Fields and Complexes}
\begin{aug}
\author{\fnms{Robert J.} \snm{Adler}\thanksref{t1,t2}\ead[label=e1]{robert@ee.technion.ac.il}
\ead[label=u1,url]{webee.technion.ac.il/people/adler}}
\author{\fnms{Omer}
  \snm{Bobrowski}\thanksref{t1,t2}\ead[label=e2]{bober@tx.technion.ac.il}
\ead[label=u2,url]{webee.technion.ac.il/people/bober/} }
\author{\fnms{Matthew S.}
  \snm{Borman}\thanksref{t2}\ead[label=e4]{borman@math.uchicago.edu}
\ead[label=u4,url]{http://www.math.uchicago.edu/$\sim$borman/} }
\author{\fnms{Eliran} \snm{Subag}\thanksref{t1,t2}
\ead[label=e3]{selirans@t2.technion.ac.il}
}
\and
\author{\fnms{Shmuel} \snm{Weinberger}\thanksref{t2,t3}
\ead[label=e5]{shmuel@math.uchicago.edu}
\ead[label=u5,url]{http://www.math.uchicago.edu/$\sim$shmuel/}}

\thankstext{t1}{Research supported in part by  US-Israel Binational
Science Foundation, 2008262.}
\thankstext{t2}{Research supported in part by NSF-SGER grant DMS-0852227.}
\thankstext{t3}{Research supported in part by NSF DMS-0805913 and by the DARPA
STOMP program.}
\runauthor{Adler, Bobrowski, Borman, Subag, Weinberger}
\affiliation{Electrical Engineering, Technion -- Israel Institute of Technology\\ Department of Mathematics, University of Chicago}
\end{aug}

\begin{abstract}
We discuss and review recent developments in the area of applied algebraic
topology, such as persistent homology and barcodes. In particular, we
discuss how these are related to understanding more about manifold learning
from random point cloud data, the algebraic structure of simplicial complexes
determined by random vertices and, in most detail, the algebraic topology  of the
excursion sets of random fields.

\end{abstract}


\end{frontmatter}

\section{Introduction}
\sec
Over the last few years there has been a very interesting and rather exciting
development in what is reputedly one of the most esoteric areas of pure mathematics:
algebraic topology. Some of the practitioners  of this subject
  are, to a considerable extent, looking out
beyond the inner beauty of their subject and seeing if they can apply it
to problems in the `real world', that is to problems outside the realm
of pure mathematics. As a result
`applied algebraic topology' is no longer an oxymoron, and although it is true that at this point sophisticated
 applications are still few and far between, there is a growing feeling
that the gap between theory and practice is closing.
We shall give more specific references below, but a very lively discussion
of this trend can be found in Rob Ghrist's review  \cite{GHRIST-BARCODES},
book in progress   \cite{ghrist2008notes} and website on a project on sensor
topology for minimal planning  \cite{Ghrist-STMP}.
 Gunnar Carlsson's webpage
\cite{Carlsson-webpage}, which describes a large  Stanford TDA (topological data
analyis) project, and a DARPA webpage \cite{DARPA-webpage}
describing a broad based project,  also help  explain the reasons why so  many people have been so attracted to this direction.

These ideas are not totally new. For example, the brain imaging community
has been using random field modelling and topological properties of these
fields for quite some time.
 For example, people like Karl Friston, a leading figure
in medical imaging, have been talking about the notion of
`topological inference' for a while (cf.\ the website \cite{Friston-website})
based in a large part on the work of the late Keith Worsley. What is new,
however, is the coordinated attack of a goodly number of high powered
mathematicians on applications.

The aim of the current paper is to describe some of the new ideas that
have arisen in applied algebraic topology and, given the interests of the
authors, exploit some of them in the setting of random fields, i.e.\ of random
processes defined over spaces of dimension greater than one. There are new
results here, albeit without proofs.  However, this paper
 is mainly review and exposition,
 with a strong bias in a particular direction, but
written in a language which we hope will be accessible to the natural
readers of this Festchrift who may (as did we until recently) find the language
of even applied topologists somewhat unfamiliar. In the final
analysis, if Larry will be happy with the final product, then
we shall be happy as well.

The paper starts in Section \ref{persistence:sec} with a discussion of one of
the central notions of applied algebraic topology, that of persistent homology
and its graphical depiction via barcodes. This is done via examples rather
than formal definitions, so it should be possible to understand the notion
of persistent homology without actually knowing what a homology group is.
(For those who do know about homology, more precise definitions
of persistent homology and barcodes are given in the appendix of Section \ref{appendix}.)
Also in Section \ref{persistence:sec} we discuss simplicial complexes as
they arise in manifold learning and also discuss the topology of
random field excursion sets.

Section \ref{fields:sec} has a brief discussion of persistence diagrams  of
excursion sets, based on simulations. (These have actually already been used
elsewhere for the analysis of brain imaging data, see \cite{Bubenik-Moo}.)
These data raise numerous  challenges for statisticians and probabilists.

Section \ref{sec:euler_integration} introduces what seems at first to be a
rather abstruse structure of  `Euler integration', but it is  very quickly
shown that not only is this a useful concept, but the key to solving a
number of quite varied problems. This is the main section of the paper.

A very brief Section \ref{complexes:sec} points out that we should have
also had more to say about random simplicial complexes, but didn't, and so
points you to appropriate references. A brief technical
appendix completes the paper in Section \ref{appendix}.

\section{Persistent homology and barcodes}
\sec
\label{persistence:sec}

In this section we are going to give a very brief and sketchy introduction to
some basic notions of algebraic topology.  A concise, yet very clear introduction to the topics that concern us can be found in \cite{bubenik2007statistical,ghrist2008notes}, while \cite{hatcher2002algebraic,vick1973homology} are good examples of a thorough coverage of homology theory. Recent excellent and quite different reviews by   Carlsson \cite{Carlsson-review,Carlsson-persistence},
Edelsbrunner and Harer \cite{Edelsbrunner-survey}, and Ghrist
\cite{GHRIST-FACTORIES,GHRIST-SENSORS,GHRIST-SENSORS2,GHRIST-BARCODES}
give a broad exposition of the basics of persistent homology.

Algebraic topology focuses on studying topology by assigning
algebraic, group theoretic, structures to topological spaces $X$.
 Thus, homology, cohomology and homotopy groups can be used to classify objects into
 classes of `similar shape'. In this paper we shall
 focus on homology. If $X$ is of dimension $N$, then it has $N+1$  homology
groups, each one of which is an abelian group. (We shall later take the coefficients
 from $\Z_2$, thereby making the groups vector spaces.) The zero-th homology $H_0(X)$ is generated by elements  that represent connected components of $X$. For $k\ge 1$ the $k$-th homology group $H_k(X)$ is generated by elements representing $k$-dimensional `loops' in $X$. The rank of $H_k(X)$, denoted by $\beta_k$, is called the $k$-th  {\it Betti number}. For $X$ compact and $k\ge 1$, $\beta_k$, measures the number of  $k$-dimensional holes in $X$, while $\beta_0$ counts the number of connected components. The Euler characteristic, a central topological quantity and homotopy invariant, is then
\beq
\label{ec-defn}
\chi (X) \ =\ \sum_{k=0}^N (-1)^k \beta_k.
\eeq

To explain the idea of persistent homology, we shall work with two examples.
The first is based on  what is known as the
`Morse filtration' of  excursion sets, the second on complexes formed from
point sets.

\subsection{Barcodes of excursion sets}
Suppose that $M$ is a nice space, that  $f:M\to\real$ is smooth, and consider the excursion, or super-level, sets
\beq
\label{Au-defn}
A_u \ \definedas \{p\in M\: f(p) \in [u,\infty)\} \
\equiv\ f^{-1}([u,\infty)).
\eeq
 Note that if $u\geq v$ then $A_{u}\subseteq A_{v}$. Going from
$u$ to $v$, components of $A_u$ may merge and new components may be born and
possibly later merge with one another or with the components of $A_u$. Similarly,
the topology of these components may change, as holes and other structures
form and disappear. Following the topology of these sets, as a function
of $u$, by following their homology, is an example of persistent homology.
The term `persistence' comes from the fact that as the level $u$ changes
there is no change in homology until reaching a level $u$ which is a
critical point of $f$; i.e.\ the topology of the excursion sets
remains static, or  `persists', between the
heights of critical points. This, of course, is the basic observation of
Morse Theory, which links critical points to homology.  However, the persistence of persistent homology goes further. For example, when two components
merge, one treats the first of these to have appeared as if it is
continuing its existence beyond the merge level.

\begin{figure}[ht!]
\begin{center}
\mbox{\scalebox{0.4}{\includegraphics*{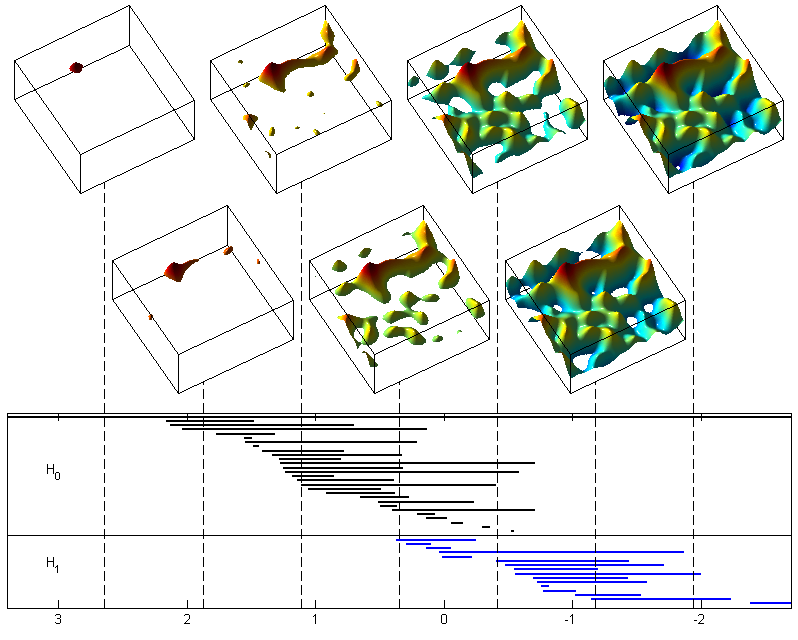}}}
\end{center}
\caption{Barcodes for the excursion sets of a function on $[0,1]^2$.
The top seven boxes show the surfaces generated by
 a 2-dimensional random field above  excursion sets $A_u$
for different levels $u$.
To determine the level for each figure,
follow the vertical line down to the scale at the bottom of the barcode.
As the vertical lines pass through the boxes labelled $H_0$ and $H_1$, the
number of intersections with bars in the  $H_0$ ($H_1$) box gives the number
of connected components (resp.\ holes) in $A_u$. Thus, at $u\sim 1.9$, $A_u$
has 4 connected components but no holes, while at $u\sim -1.2$, $A_u$ has only
1 connected component, but 9 holes. The horizontal lengths of the bars
indicate how long the different topological structures (generators of the
homology groups) persist. Computation of the barcodes was carried out in
Matlab using  Plex (Persistent Homology Computations) from Stanford
 \cite{PLEX}.}
\label{fig:2d}
\end{figure}

A useful way to describe persistent homology is via the notion of barcodes.
Assuming that $\dim (M)=N$, we also have, from the smoothness of $f$,
that, if $A_u$ is non-empty, then  $\dim (A_u)$ will typically also be $N$.  A barcode for the excursion sets of $f$ is then a collection of $N+1$ graphs,
one for each collection of homology groups of common order. A bar in the
$k$-th graph, starting at $u_1$ and
ending at $u_2$ ($u_1\geq u_2$) indicates the existence of a generator
of $H_k(A_u)$ that appeared at level $u_1$ and disappeared at level $u_2$.
An example is given in Figure \ref{fig:2d}, in which the function $f$ is
actually the realisation of a smooth random field on the unit square, an
example to which we shall return  later.

\begin{figure}[ht!]
\begin{center}
\mbox{\scalebox{0.4}{\includegraphics*{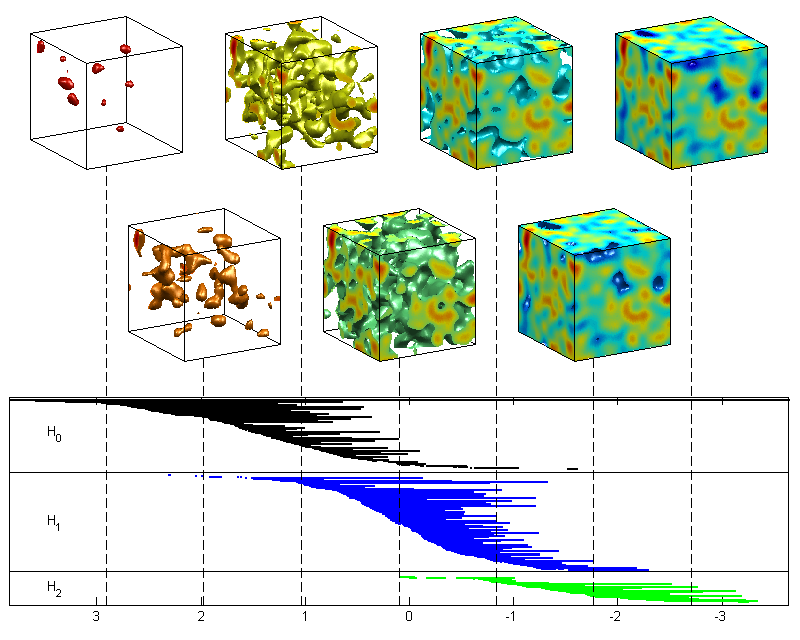}}}
\end{center}
\caption{Barcodes for the excursion sets of a 3-dimensional random field.
The barcode diagram is to be read as for Figure \ref{fig:2d},
with two differences: The top 7 boxes now  display the excursion
sets themselves and the values of the field are colour coded. Furthermore,
there are now three homology-groups/barcode-boxes, representing connected
components, handles, and holes.}
\label{fig:3d}
\end{figure}

Figure \ref{fig:3d} is even more impressive, since it shows a three
dimensional example.  Note that, as opposed to the 2-dimensional case,
it is almost impossible to
say anything about topology just by looking at the boxes with the excursion sets
at the top of the figure, but there is a lot of immediate visual information
available in the barcodes. This phenomenon becomes even more marked as the
dimension $N$ of the parameter space increases. While it
 may be impossible to imagine what a five dimensional excursion set looks
like, it is easy to look at a barcode with six sets of bars for the six
persistent homologies.

\subsection{Point sets and manifold learning}
\label{pointsets:subsec}
\begin{figure}[ht!]
\begin{center}
{\scalebox{0.25}{\includegraphics*{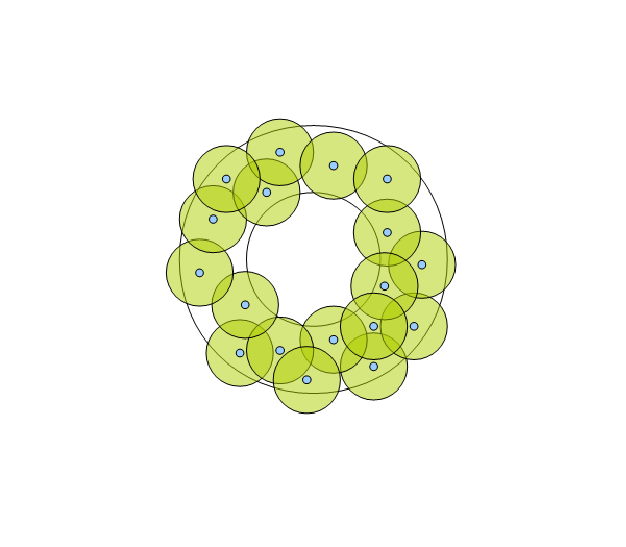}}}
{\scalebox{0.25}{\includegraphics*{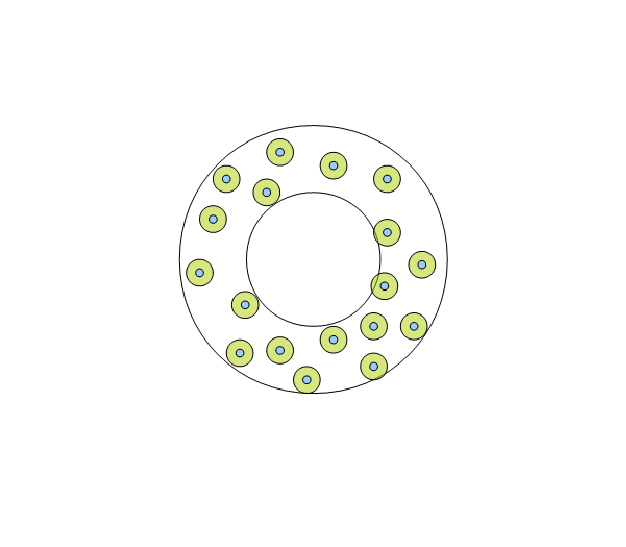}}}
{\scalebox{0.25}{\includegraphics*{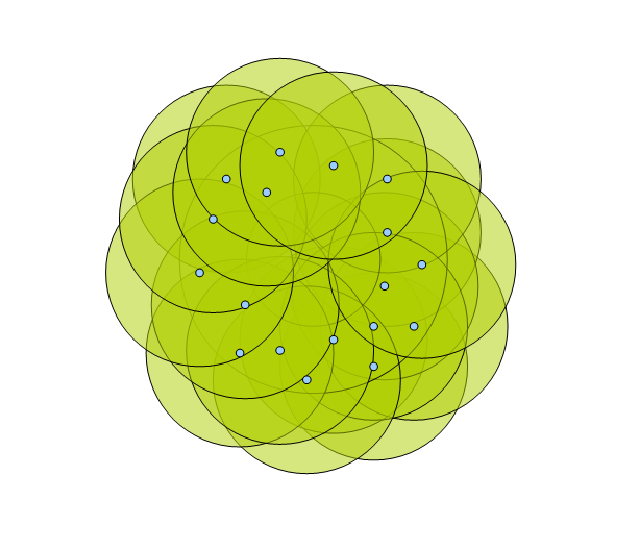}}}
\end{center}
\caption{Trying to capture the homology of an annulus (where $\beta_0 = 1,\, \beta_1 = 1$) from a union of balls of various radii around a random sample of points from the
annulus.
 A good choice of radius recovers the correct homology in the first case.
If the radius chosen is too small,  the union of balls has the same homology as $n$ distinct points ($\beta_0 = n, \beta_1 = 0$). If the radius chosen is too big, the union is contractible ($\beta_0 = 1, \beta_1 = 0$).}
\label{balls-fig}
\end{figure}

Consider the following situation. Let $X$ be an unknown subset of $\R^d$
 with  finite Lebesgue measure and let $X_1,\ldots,X_n$ be $n$ independent random samples uniformly distributed on $X$. We would like to study the homology of $X$ using only these  random points. When $X$ is a manifold, this is typically
referred to as {\it manifold learning}. In many cases we can find an $\e$ for which the union of balls
\beq
\label{unionofballs}
U\ =\ \bigcup_{i=1}^n{B_\e(X_i)}
\eeq
 is homotopy equivalent to $X$ (and hence has the same homology). However, we do not know,  a priori,  what is the correct choice of $\e$. An example is given in Figure  \ref{balls-fig},
in which $X$ is a two-dimensional annulus. If $\e$ is chosen to be too small  then $U$ is homotopy equivalent to the union of $n$ distinct points (and hence contains no information on $X$). On the other hand, choosing $\e$ to be too big gives us a $U$ that is a large, contractible blob, which again tells us nothing about $X$. But, as with Goldilock's porridge, choosing $\e$ `just right',
recovers an object topologically equivalent to the annulus.
Persistent homology overcomes this sensitivity to the choice of $\e$
by considering a  range of possible values of $\e$, much as we did
with  the levels of excursion sets in the previous example, but with the
aim of learning about the topology of $X$ from the barcodes.
  The key assumption is that homology elements that `live longer' (or, \textit{persist}) are more likely to represent homology elements of $X$, whereas the
shorter ones  are just `noise'.

To describe this in a little more detail we need the notion of simplicial
complexes.

\subsection{Simplicial complexes}
\label{simpcomps:subsec}
We are not going to give a definition of
simplicial complexes here, but rather shall describe  two classic ways to  construct abstract simplicial complexes from a given set of points in a metric space.

\begin{definition}[The \cech Complex]\label{def:cech_complex}
Let $\cP = \set{x_1,x_2,\ldots}$ be a collection of points in a metric space $X$. Construct an abstract simplicial complex $C(\cP, \e)$ in the following way:
\begin{enumerate}
\item The $0$-simplices are the points in $\cP$,
\item An $n$-simplex $[x_{i_0},\ldots,x_{i_n}]$ is in $C(\cP,\e)$ if $\bigcap_{k=0}^{n} {B_{\e}(x_{i_k})} \ne  \emptyset$,
\end{enumerate}
where $B_\e(x)$ is the ball of radius $\e$ around $x$.
The complex $C(\cP,\e)$ is called the \cech complex attached to $\cP$ and $\e$.
\end{definition}

\begin{definition}[The Vietoris-Rips Complex]\label{def:rips_complex}
Let $\cP = \set{x_1,x_2,\ldots}$ a collection of points in a metric space $X$. Construct an abstract simplicial complex $R(\cP, \e)$ in the following way:
\begin{enumerate}
\item The $0$-simplices are the points in $\cP$.
\item An $n$-simplex $[x_{i_0},\ldots,x_{i_n}]$ is in $R(\cP,\e)$ if $B_{\e}(x_{i_k}) \cap B_{\e}(x_{i_m}) \ne \emptyset$ for every $0\le k < m \le n$.
\end{enumerate}
The complex $R(\cP,\e)$ is called the Rips complex attached to $\cP$ and $\e$.
\end{definition}

From these definitions it is obvious that $C(\cP,\e) \subset R(\cP,\e)$. In addition, it is proved in \cite{GHRIST-SENSORS} that $R(\cP,\e') \subset C(\cP,\e)$ for $\e / \e' \ge \sqrt{2d/(d+1)}$. In other words, a \cech complex can be `approximated' by Rips complexes. This fact is used in computational applications, since working with Rips complexes is much more efficient than with \cech complexes.

There are occasions when Rips and \cech complexes coincide, as is
the case when $X$ is Euclidean but the metric is the $L^\infty$ rather than
the more standard $L^2$ norm.
In many statistical applications the choice of metric on $X$ may be dictated
by optimality considerations rather than `natural' geometry.

The main importance of the \cech complex and its relevance to homology theory, is given in the next theorem.

\begin{Theorem}[The Nerve Theorem]\label{thm:nerve}
Suppose that the intersections $\bigcap_{x\in \cal P'} B_{\e}(x)$ are either empty
or contractible for any subset $\cP'$ of $\cP$. Then the \cech complex $C(\cP,\e)$ is homotopy equivalent to $\bigcup_{x\in \cP} B_{\e}(x)$. In particular,
if $X$ is a finite dimensional normed linear space, or a compact Riemannian
manifold with convexity radius greater than $\e$,
and  if $\set{B_{\e}(x)}_{x\in \cP}$ is a cover of the space $X$,
then $C(\cP,{\e})$ is homotopy equivalent to $X$.
\end{Theorem}

The main consequence of the Nerve Theorem is that in order to study the
homology of the {\it topological} space $\bigcup_{x\in \cP} B_{\e}(x)$, we can study the homology of the {\it combinatorial} space $C(\cP,{\e})$.
 This fact can be useful in proving theoretical results, but its main contribution is to computational applications.

\begin{figure}[ht!]
\begin{center}
\mbox{\scalebox{0.5}{\includegraphics*{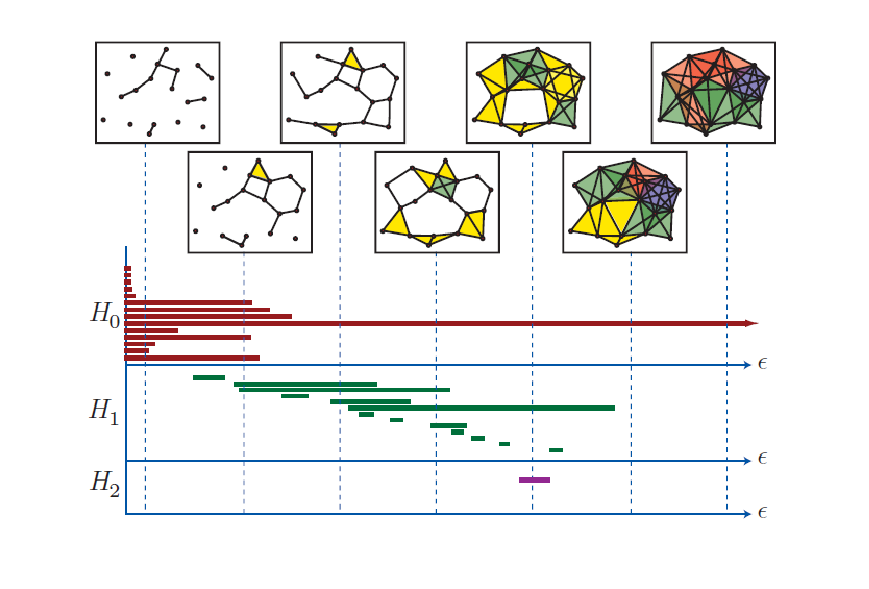}}}
\end{center}
\caption{The barcode of a Rips complex, taken from \cite{GHRIST-BARCODES}. The points were sampled from an annulus in $\R^2$. We see that there is a single $H_0$ bar that persists forever. This bar represents the single connected component of the annulus. In $H_1$ we see a couple of dominant bars indicating that the sample space contains holes. The longest bar actually represents the real hole of the annulus. In $H_2$ there is nothing significant and indeed $\beta_2=0$ in this case.}
\label{fig:barcode}
\end{figure}

With these definitions behind us, Figure \ref{fig:barcode} gives a nice example of how barcodes describe the topology of an annulus
($\beta_0=1,\ \beta_1=1,\beta_2=0$) in $\real^2$, when 17 points are sampled from it and Rips complexes are computed for a range of $\e$.

\section{Random field simulations}
\sec
\label{fields:sec}

In this section we want to consider the persistent homology of random
field excursion sets. In particular, we would like to understand something
about the distributional properties of their barcodes.

The random fields behind the barcodes of Figures \ref{fig:2d} and
\ref{fig:3d} were taken to be mean zero, Gaussian, over the parameter set
$[0,1]^2$ and with covariance function $R(p)=\exp(-\alpha\|p\|^2)$. This is
a stationary, isotropic, and infinitely differentiable random field, and the
starting test case for all theories. We took $\alpha=100$.

\begin{figure}[ht!]
\begin{center}
{\scalebox{0.3}{\includegraphics*{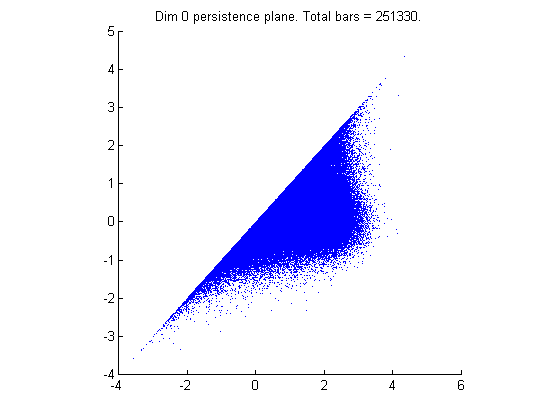}}}
{\scalebox{0.3}{\includegraphics*{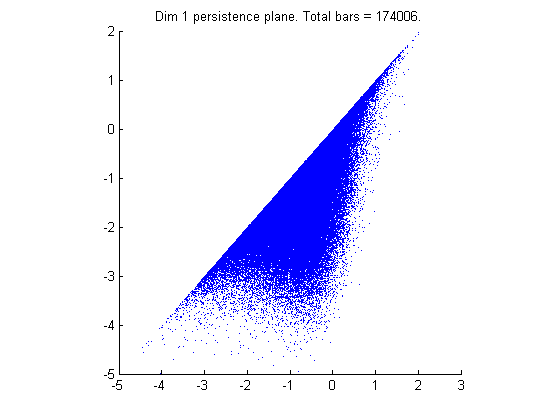}}}
\end{center}
\caption{Persistence diagrams for 10,000 simulations of an isotropic
random field on the unit square. Note that the diagrams for $H_0$ and
$H_1$ seem quite different.}
\label{fig:persistence}
\end{figure}

We ran 10,000 simulations of this field, calculating 10,000 barcodes. In order to represent the data in a reasonable fashion, we used {\it persistence
diagrams} rather than barcodes. To form a persistence diagram from the
bars in $H_k$, one simply replaces each bar by a pair $(x,y)$, where $x$ is
the level at which the bar begins and $y$ the level at which it ends. Thus
$x>y$ and the pair $(x,y)$ lies in a half plane. In Figure
\ref{fig:persistence} the corresponding persistence diagrams for the
complete simulation data are shown for $H_0$ and $H_1$.

\begin{figure}[ht!]
\begin{center}
\resizebox{2.3in}{3.0in}{\includegraphics*{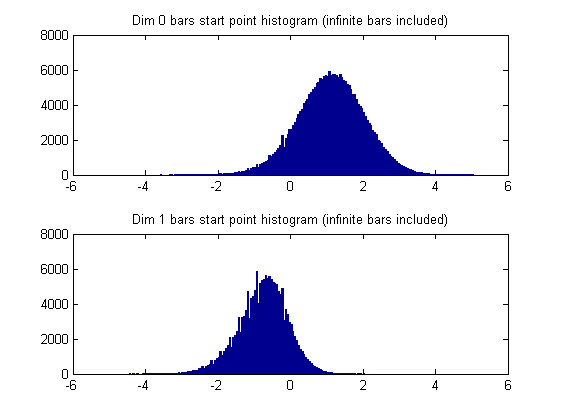}}
\resizebox{2.3in}{3.0in}{\includegraphics*{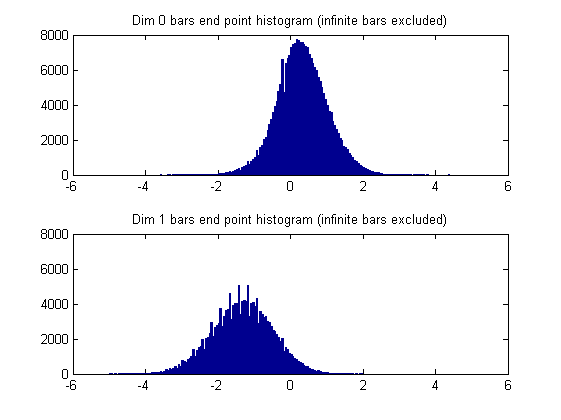}}
\end{center}
\caption{Empirical distributions of start and end points of bars for the
Gaussian field of Figure \ref{fig:persistence}.}
\label{fig:end-start}
\end{figure}

Additional information on the barcodes is given in Figure \ref{fig:end-start}.
What is shown there are the (marginal) distributions of the start and end
 points of the barcodes for $H_0$ and $H_1$ from the same simulation. A simple
application of Morse theory, or, in this simple two dimensional setting,
a little thought,  leads to the realisation that the start points of the $H_0$ bars
are all  heights of local maxima of the field, while  the end points of the  $H_1$ bars correspond to local minima. These distributions have been well studied (although their precise form is not known) in the general theory of
Gaussian random fields. The remaining start and end
points correspond to different types of saddle points of the random field.
However, what differentiates between the end point of a $H_0$ bar and
the start point of a $H_1$ bar is global geometry and is not determined by
the local behaviour of the field.

It would be interesting to know more about the real distributions lying
behind Figures  \ref{fig:persistence}  and  \ref{fig:end-start}, but at
this point we know very little. An interesting aspect is the asymmetries between
the start points of the $H_0$ bars (local maxima) and the end points of the $H_1$
bars (local minima, as well as between the two sets of saddle points. We imagine
 that this is due to boundary effects and would disappear if the simulation
had been carried out on a closed manifold.

There are some things that we do know, however, and we turn to them next.
Firstly, however, we need to make a small digression.

\section{Euler Integration} \label{sec:euler_integration}
\sec

Before we introduce the Euler integral, we need to define the Euler characteristic for noncompact spaces. For a compact space $X$, we already defined the Euler
characteristic $\chi(X)$  in \eqref{ec-defn} as an alternating sum of
Betti numbers; viz.\ as an alternating sum of ranks of the homology groups
 $H_k(X)$.

In this setting the Euler characteristic is a homotopy invariant and is additive in the sense that
\begin{equation}\label{eq:ec_additivity}
	\chi(A \cup  B) = \chi(A) + \chi(B) - \chi(A \cap B).
\end{equation}

In extending the definition of the Euler characteristic to noncompact spaces, if one uses the definition
of the  Euler characteristic as the alternating sum of $\rk H_k(X)$, then additivity is lost
(consider $[0,1] = [0,1) \cup \{1\}$).  Therefore the definition of the Euler characteristic we shall
use for noncompact spaces is
\[
	\chi(X) = \sum_k (-1)^k \rk H_k^{\lf}(X) = \sum_k (-1)^k (\#\mbox{ of open $k$-simplices in X}),
\]
where $H_k^{\lf}(X)$ is called the locally finite homology (see \cite[Chapter 3]{hughes1996}).  Since
$H_k^{\lf}(X) = H_k(X)$ for compact spaces, this extends the definition of the {\EC} to noncompact spaces
such that additivity is preserved, but we lose homotopy invariance, although
it is still a homeomorphism invariant.
This {\EC} can also be computed by decomposing the space into a union of open $k$-simplices and points.
For example, $\chi((0,1)) = -1$, $\chi([0,1)) = 0$ and $\chi([0,1]) = 1$.

\subsection{The Euler Integral}
Since the {\EC} is an additive operator on sets (cf.\ \eqref{eq:ec_additivity}), it is tempting to consider $\chi$ as a measure and integrate with respect to it. The main problem in doing so is that $\chi$ is only finitely additive.

At first (cf.\ \cite{viro1988some}), integration with respect to the {\EC} was defined for a small set of functions called \textit{constructible functions} defined by
\[
C\!F(X) = \set{\left.h(x) = \sum_{k=1}^n a_k \indic_{A_k}(x)\,\right|\, n\in \mathbb{N} ,a_k\in \Z,\, A_k \textrm{ is tame}},
\]
where `tame' means having a finite Euler characteristic.
For this set of functions we can define the {\EI} similarly to the Lebesgue integral. Let $h(x) = \sum_{k=1}^n a_k \indic_{A_k}(x)$ and define
\[
\int_X h d\chi\  \triangleq\  \sum_{k=1}^n a_k \chi(A_k).
\]
This integral has many nice properties, similarly to those of the Lebesgue integral, such as linearity and a version of the Fubini theorem (cf. \cite{ghrist2008notes,viro1988some}).
However, due to the lack of countable additivity one cannot easily
continue from here by approximation to integrate more general functions.

Nevertheless, in \cite{baryshnikov2009euler} two possible extensions were suggested for the {\EI} of real valued functions. We shall not go into details of
the constructions here, but rather use one of the properties of these
extensions to define what in  \cite{baryshnikov2009euler} was called an upper Euler integral and which we,
for simplicity,
shall call an Euler integral. Thus, we  define the Euler integral by
\beq
	\int_X f \, d\chi \ \definedas \
  \int_{u=0}^\infty\sbrk{ \chi(f > u) - \chi(f \le -u)}\, du,
\label{Eulerint-defn}
\eeq
where $\chi(f > u) \triangleq \chi\param{f^{-1}(u,\infty)}$ and
$\chi(f \leq -u) \triangleq \chi\param{f^{-1}[-u,\infty)}$.
These integrals are defined for what are known as `tame' functions.
See \cite{baryshnikov2009euler} for details.

Unfortunately, these extensions of the {\EI} have many flaws, of which the most prominent one is the lack of additivity. For example, a simple computation shows that for $X=[0,1]$
\[
\int_X x \,d\chi + \int_X (1-x) \, d\chi = 1+1 = 2 \ne 1 = \int_X 1 \, d\chi .
\]
Nevertheless, these integrals still have interesting properties and some intereting applications. Here is one.

\subsection{An application of the Euler integral} \label{sec:target_enum}
This application was suggested in \cite{baryshnikov2007target}.
Suppose that an unknown number of targets are located in a space $X$ and each target $\alpha$ is represented by its support $U_{\alpha}\subset X$. Suppose also that the space $X$ is covered with sensors, each reporting only the number of targets it can sense, but with no ability to distinguish between targets.
 Let $h:X \to \Z$ be the \textit{sensor field}, i.e.\! \[h(x) =\#\set{\textrm{targets activating the sensor located at $x$}}.\] The following theorem states how to combine the readings from all the sensors and get the exact number of targets.

\begin{Theorem}[Baryshnikov and Ghrist, \cite{baryshnikov2007target}]
\label{thm:target_count}
If all the target supports $U_{\alpha}$ satisfy $\chi(U_{\alpha}) = \gamma$ for some $\gamma\ne 0$, then
\[
\#\set{targets} = \frac{1}{\gamma}\int_X h d\chi.
\]
\end{Theorem}

Note that we do not need to assume anything about the targets other than they all have the same nonzero {\EC}. For example, we need not assume that they are all convex or even have the same number of connected components. On the other hand, the theorem assumes an ideal sensor field, in the sense that the entire (generally continuous) space $X$ is covered with sensors which register only what happens at the point at which they are placed. In \cite{baryshnikov2009euler} more realizable models using the upper and lower Euler integrals are discussed.

Assume now that the readings from the sensors  are contaminated by a Gaussian (or Gaussian related) noise $f(x)$. Under these conditions it can be proved that
\[
\int_X (h+f)\,d\chi = \int_X h\,d\chi + \int_X f \,d\chi.
\]
Denoting $s= \int_X h \,d\chi$ (deterministic signal), $n=\int_X f \,d\chi$ (noise) and $y=\int_X (h+f) \,d\chi$ (measurement), this is a classic signal plus
 noise problem  (i.e.\! $y = s+n$). In particular, in order to estimate $s$ from $y$,
it would be nice, in view of Theorem \ref{thm:target_count}, to be able to compute some distributional properties of the Euler integral of a Gaussian random field. We shall limit ourselves to  computing the expectation and shall  turn to this
after a
few words on the Gaussian kinematic formula.

\subsection{The Gaussian kinematic formula}
Suppose that $M$ is an $N$-dimensional, $C^2$, Whitney
 stratified manifold satisfying some mild side conditions (cf.\ \cite{RFG}
for details) and $D$  a similarly nice  stratified submanifold of $\Rk$. Let
 $f=(f^1,\dots,f^k):M\to\real^k$ be
 a vector valued random process, the components of which are independent,
 identically distributed,  real valued, $C^2$, centered,
unit variance, Gaussian processes.
Using $f$, define a  Riemannian  metric on $M$ by setting
\beq
\label{metric}
g_x(X,Y)\ \definedas\ \E\{(Xf^i_x)\,( Yf^i_x)\},
\eeq
for any $i$ and for $X,Y\in T_xM$, the tangent space to $M$ at $x\in M$,
and use this to define the \LK\ curvatures,  $\lips_j$, $j=0,\dots,N$ on
$M$.  For example, if $M$ is a manifold without boundary, then these are
given by
\begin{equation}
\label{Lipschitz:Killing:Without}
\lips_j(M) \ =\
\frac{1}{(2 \pi)^{(N-j)/2} ((N-j)/2)!}\int_M
\Tr^M(-R)^{(N-j)/2}\, {\rm Vol}_g,
\end{equation}
when $N-j \geq 0$ is even, and 0 otherwise.
Here ${\rm Vol}_g$ is the volume form of the Riemannian manifold $(M,g)$,
$R$ is the curvature tensor and
 $\Tr^M$ the trace operator on the algebra of double forms on $M$.
For simple Euclidean spaces, with various orderings and normalisations,
the \LK\ curvatures are also known as  Quermassintegrales, Minkowski or Steiner functionals, integral curvatures, and  intrinsic volumes. Note that
$\lips_N(M)\equiv {\rm Vol}_g(M) $ is
the Riemannian volume of $M$ and $\lips_0(M)\equiv \chi(M)$
is its Euler characteristic.

The Gaussian kinematic formula (hereafter GKF) was due
 originally to Taylor in \cite{Taylor-thesis} (but
for  the  form below see \cite{RFG,TAannals})  and states that
\beq
\label{p:main:equation}
\E\left\{\lips_i\left(M\cap f^{-1}(D)\right)\right\}\ =\
\sum_{j=0}^{N -i} \sqbinom{i+j}{j} (2\pi)^{-j/2}\lips_{i+j}(M)
\Min_j^{\gamma}(D).
\eeq
The combinatorial
coefficients here  are the standard `flag coefficients' of integral geometry,
given by
 \beqq
\sqbinom{n}{j}\  =\   \binom{n}{j} \frac{\omega_n}{\omega_{n-j} \; \omega_j},
\eeqq
where $\omega_n$ is the volume of the
unit ball in $\real^n$. The $\Min_j^\gamma(D)$, known as  the  Gaussian Minkowski functionals of $D$, are determined via the tube expansion
\beq
\label{tubes:gaussexp1:equn1}
\P\left\{ f(x)\in\{y\:d(y,D)\leq \rho\}\right\} \ =\  \sum_{j=0}^{\infty}
\frac{\rho^j}{j!} \Min_j^\gamma(D),
\eeq
where $x$ is any point in $M$ and $d$ is the usual Euclidean distance from a
point to a set.

One could devote a book  to this formula and, indeed, such a book exists.
So we shall refer you to  \cite{RFG} for all needed technical details.

We note only one pertinent fact, for immediate use.  Taking $j=0$ in
\eqref{p:main:equation} gives the expected Euler characteristic of excursion
sets as a simple, closed form expression that can be readily calculated in
many interesting cases. Again, see \cite{RFG} for details.

\subsection{The Euler integral of a Gaussian random field}
 \label{sec:ei_grf}
Returning to the signal plus noise problem of Section
\ref{sec:target_enum}, we  can  formulate the first step towards  its solution.

Let $M$ be a nice, tame, space. (The definition of `tame' can be found in
\cite{RFG}.) Let $f$ be a  random field.
Here is a striking result, due to Bobrowski and Borman \cite{BandB1}:

\begin{Theorem}\label{thm:mean_ei_general}
Let $M$ be an $N$-dimensional tame stratified space and let $f:M\to\R^k$ be a $k$-dimensional Gaussian random field satisfying the GKF conditions. Let $G:\R^k\to\R$ be piecewise $C^2$ and let $g=G\circ f$. Setting $D_u = G^{-1}(-\infty,u]$ and assuming that $\abs{\int_{\R}{\MM_j(D_u) du}}<\infty$,  we have
\begin{equation}\label{eq:mean_ei}
\mean{\int_M{g \,d\chi}} = \chi(M)\mean{g} - \sum_{j=1}^N{(2\pi)^{-j/2}\LL_j(M)\int_{\R}{\MM_j^\gamma(D_u) du}},
\end{equation}
where $\mean{g} := \mean{g(t)}$. ($g(t)$ has constant mean).
\end{Theorem}

While, on the one hand, this is not a difficult result to prove, given
the GKF and \eqref{Eulerint-defn}, it was completely unforseen until
discovered and has a number of interesting and potentially deep implications.

The main difficulty in applying Theorem \ref{thm:mean_ei_general}
 lies in computing the Minkowski functionals $\MM_j^\gamma (D_u)$.
A simple example is given in the following case:

\begin{Theorem}
\label{thm:mean_ei_real}
Let $M$ be an $N$-dimensional tame stratified space and let $f:M\to\R$ be a real valued Gaussian random field satisfying the GKF conditions. Let $G:\R\to\R$ be  piecewise $C^2$ and let $g=G \circ f$. Then
\[
\mean{\int_M{g \,d\chi }} = \chi(M)\mean{g} + \sum_{j=1}^N{(-1)^j\LL_j(M)} \frac{\iprod{H_{j-1}, (\sign(G'))^j G'}}{(2\pi)^{j/2}}.
\]
\end{Theorem}

In the theorem, for $n\ge 0$ the $n$th Hermite polynomial $H_n $ is defined by
\[
	H_n(x) = (-1)^n \varphi^{-1}(x) \frac{d^n}{dx^n} \varphi(x),
\]
where
$\varphi$ is the standard Gaussian  density, $e^{-x^2/2}/\sqrt{2\pi}$. The inner product is given by
\begin{equation}\label{eq:inner_prod}
	\iprod{f,g} = \int_{\R}{f(x)g(x)\vp(x)}dx,
\end{equation}
and we use the convention, required below, that
\[
	H_{-1}(x) = \varphi^{-1}(x)\int_x^\infty \varphi(u) du.
\]

In the case that the function $G$ is strictly monotone, we have an even simpler form.

\begin{Corollary}
\label{cor:mean_ei_real_monotone}
Let $f$ be as in Theorem \ref{thm:mean_ei_real} and $G$ be a strictly increasing function, then
\[
\mean{\int_M{g \,d\chi }} = \sum_{j=0}^N {(-1)^j\LL_j(M) \frac{\iprod{H_j,G}}{(2\pi)^{j/2}}}.
\]

If $G$ is strictly decreasing then,
\[
\mean{\int_M{g \,d\chi }} = \sum_{j=0}^N {\LL_j(M) \frac{\iprod{H_j,G}}{(2\pi)^{j/2}}}.
\]
\end{Corollary}
Finally, taking $G$ to be the identity function yields

\begin{Corollary}\label{cor:expected_ei_grf}
Let $f$ be as in Theorem \ref{thm:mean_ei_real}, then
\[
\mean{\int_M{f \,d\chi }} = -\frac{\LL_1(M)}{\sqrt{2 \pi}}.
\]
\end{Corollary}

Further details and further examples, including results for $\chi^2$ and $F$ random
fields,  can be found in \cite{BandB1}.

\subsection{Persistent homology of Gaussian excursion sets}
 \label{sec:gaussian_fields}

We now return to excursion sets, which would seem to have been forgotten
in all the discussion on Euler integrals. It turns out that this was not the
case, but that Euler integrals and Theorem \ref{thm:mean_ei_real} contain a
lot of information on the persistent homology of Gaussian excursion sets.

First, however, some notation and a definition:
Suppose we have a barcode, which we shall denote by $\cal B$.
 Denote the
individual bars in $\cal B$ by $b$, their lengths by $\ell (b)$, and the
degree  of the homology group to which belongs the generator that they
represent by $\mu (b)$.
\begin{definition}\label{def:ec_barcode}
	The \emph{Euler characteristic} of a barcode $\cal B$ with no bars of infinite length is
	\[
		\chi(\cal B) \ \definedas \
 \sum_{b \in \cal B} (-1)^{\mu(b)} \ell(b).
	\]
\end{definition}

It turns out that this topological quantity can actually be written in terms
of Euler integrals. For convenience, and adopting the prejudices
 of topologists rather than probabilists and statisticians, we shall
consider the barcodes of incursion rather than excursion sets, or,
equivalently, sub-level sets rather than super-level sets. That is, we
replace excursion sets  of \eqref{Au-defn} by
\beq
\label{Atildeu-defn}
\widetilde A_u \ \definedas \{p\in M\: f(p) \in (-\infty, u]\} \
\equiv\ f^{-1}((-\infty, u]).
\eeq
Then \cite{BandB1}  showed that, if
 $\cal B (f,u)$ denotes the barcode of  $\widetilde A_u$ for  tame $f$,
\beq
\label{barsandeuler}
 \chi(\cal B(f,f_{\max})) = f_{\max}\chi(M)- \int_M f \,d\chi .
\eeq

Combining this with Theorem \ref{thm:mean_ei_general} yields

\begin{Theorem}
\label{expected-barcodes}
Let $f:M\to \R^k$ be a Gaussian random field satisfying the GKF conditions, $G\in C^2(\R^k,\R)$ and $g=G\circ f$. Then
\beqq
\mean{\chi(\cal B(g,g_{\max}))}  =  \chi(M)\param{\mean{g_{\max}}-\mean{g}}+ \sum_{j=1}^N{(2\pi)^{-j/2}\LL_j(M)\int_{\R}{\MM_j(D_u) du}}.
\eeqq
If $f$ is real, then
\beqq
\mean{\chi(\cal B(f,a))} =  \chi(M)\param{\varphi(a) + a \Phi(a)} + \varphi(a) \sum_{j=1}^N{(2\pi)^{-j/2}\LL_j(M) H_{j-2}(-a)},
\eeqq
for any $a$.
\end{Theorem}

The way we have presented this result, as a `natural' consequence of a
reasonably `straightforward'  Theorem \ref{thm:mean_ei_general}, substantially
underplays its importance and novelty. It has a number of interesting
corollaries, for which we send you to the original paper \cite{BandB1}. But
its main contribution lies in its very existence, connecting, as it does,
between probabilistic objects and their homological structure.

As one of our colleagues/teachers recently stated: ``I can think of no two
topics in mathematics further away from one another than probability and
algebraic topology. There is probably no way to connect them.''
Yet here, in Theorem \ref{expected-barcodes}, is an elegant connection, one
of the first of its kind.

\section{Random Geometric Complexes}
\label{complexes:sec}
\sec
In Sections \ref{pointsets:subsec} and \ref{simpcomps:subsec} we motivated
the idea of persistent homology and barcodes with examples
from manifold learning and random simplicial complexes. Despite this, we
shall not go into  detail, but shall rather describe some general issues and
give you a few  useful references to this area, which is also
currently undergoing rapid development.

\subsection{Manifold learning}

We already mentioned manifold learning briefly in Section
\ref{pointsets:subsec} via the example of trying to identify an annulus, or at
least its homology, from  a simplicial complex built over the sample points.
The subject of manifold learning goes, obviously, well beyond such an example,
and examples of algorithms for `estimating' an underlying manifold from
a finite sample abound in the statistics and computer science literatures.
Very few of them, however, take an algebraic point of view, which is what
we have stressed in this paper.

One  contribution in the spirit of this paper is  \cite{niyogi2008finding} by
Niyogi, Smale, and Weinberger, who studied the problem of estimating
smooth manifolds from finite samples.
They showed that  in sampling from a high dimensional
manifold, if the sampling is dense enough then, with high probability,
the set \eqref{unionofballs} deformation retracts to the manifold and so has
the same homology.  This implies that long
persistence intervals, once one has enough sample points, are very
likely to correctly compute the homology of a submanifold.

Of course, one of the most important issues in dealing with data is
noise. In the setting of manifold learning  this translates to the sample
points possibly not coming from  the submanifold that
theoretically models the phenomenon because of experimental, measurement,
or other  error.
In \cite{NSW2} the same authors treat this issue, as does
\cite{CC-SM} from a different and enlightening point of view.

In a complementary but related direction \cite{BCKL,bubenik2007statistical}
apply persistence techniques to the nonparametric study of functions on a
given manifold.

\subsection{Random complexes}
We now return to the  \cech and Rips complexes of Section \ref{simpcomps:subsec}.

To get a feeling for the phenomena that occur as we approximate a
manifold by the union of balls, it is perhaps enlightening to consider
the situation, for a fixed $\e$, of the evolution of the homology
of the \cech complexes as the number of points, $n$, grows.
The points themselves we assume are chosen uniformly, at random,
 on the manifold.

  For $n$ small, the
balls do not intersect, but as $n$ grows intersections begin to occur and
small finite graphs appear.  Assuming that $\e$ is sufficiently
small that  all $\e$-balls have about the same
volume, it is easy to compute the expected number of times a
particular graph arises.  This leads to complicated integrals, but
investigating them leads to the belief
that   $k$-homology (for a \cech complex) is most likely to occur as a result of
the occurrence of boundaries of $(k+1)$-simplices. That is, it
 requires $k+1$ points to be
close to each other
(at scale $\e$) but not to fill in.  Aside from  a constant factor, the
probability that this happens  should be $(\e^d)^k{n \choose k+1}$.  For this
probability to be non-trivial, one requires that $n$ is
$O(1/\e^{dk/(k+1)})$. In other words, it is for $n$ of this size that one
begins to see interesting $k$-homology.

As $n$ continues to grow, there will be a point
where the data covers a nontrivial percentage of the volume, a point at which
phenomena related to percolation occur. Finally, there is a
reversal when the $\e$-balls fill almost all of the manifold,  all
 extra homology dies,  and ultimately we obtain  the correct
calculation of homology.

To get a feeling for the end game, it is worthwhile to compute the
expected Euler characteristic of the union of $n$  $\e$-balls in, say
a flat torus.  For simplicity, as we are only giving a heuristic,
let's avoid the complications of Euclidean metrics and consider an
$L^\infty$ metric on the torus.  In that case, a straightforward
inclusion exclusion argument (see \cite{Okun} for this in a Poisson model in
Euclidean space and \cite{Mecke-Wagner}
 for the use of kinematic formulae to obtain
the relevant formulae in the case of genuine round balls), together with a
generating function argument, give the formula for
$\E(\chi_{ n, \e, d})$, where $\chi_{ n, \e, d}$ is the Euler characteristic of the union of $n$ $\e$-balls and $d$ is the dimension of the torus, as follows:

 Let $\tau = (2\e)^d$.  Then
\beqq
\E\left\{\chi_{ n, \e, d}\right\}
\ = \ \begin{cases}
 n(1-(\tau)^{n-1})  &d=1\\
\frac{d}{d\tau} \left(\tau \E\left\{\chi_{ n, \tau, d-1}\right\}\right)
& d\geq 2.
\end{cases}
\eeqq
  One thus obtains that the
Euler characteristic is approximately 0 (for the last time and so implying
coverage) when $n$ is around $(1/\tau)\log(1/\tau)$ plus lower order terms.  See
\cite{Flatto-Newman} for more details.  Interestingly, it follows from the work of
\cite{penrose2003random} that the phase transition for the giant component to
 form,  in the sense
of random graph theory,  is (asymptotically)
at $2^{-d}$ times this number. That is, the computation of components seems to be
correct.  It also seems extremely
likely that there are phase transitions at other multiples of this
fundamental scale where the other homology groups are correctly
computed.

Many more details of the phenomena for small $n$ and the percolation
range, the relevant central limit theorems for homology and some
valuable information about persistence intervals in Rips and \cech
situations can be found in recent papers of Kahle  \cite{kahle2009random} and
Bobrowski and Borman \cite{BandB2}.  They
combine probabilistic tools with Morse theory to give rigorous proofs
of these phenomena.  We mention also an early lecture by Diaconis,
available on the web \cite{Diaconis-web}, which
  suggests the general outline of this picture,
and a forthcoming paper \cite{BW-inprep}  that also deals with some
aspects of this
problem in a metric-measure setting (relevant to situations where the
distribution of points in nature does not follow the Riemannian
uniform measure).

\section{Technical appendix}
\label{appendix}

\sec
As promised, we shall now be a little more formal and explain what persistent
homology really is. For this, however, we shall need to
assume that the reader has a
basic working familiarity with the theory of simplicial homology.
We shall also, for simplicity, take all homologies over the group $\Z_2$.
This is the material of  Section \ref{persistenth:subsec}.

The second subsection of the appendix explains how to carry these definitions
over to random complexes. To complete the Appendix we should have really
added a section on how one turns the excursion sets of a continuous random
field into a random filtered complex. This, as you will be able to guess,
after reading the first two subsections, is done by discretizing the
 parameter space of the random field and then thresholding the field at
various levels in order to obtain the simplices of the filtered complex.
You can find details of this in the report \cite{Subag1}.

\subsection{Persistent Homology}
\label{persistenth:subsec}
We start by considering growing sequences
of simplicial complexes that grow in the following manner.
\begin{definition}
\label{filtered complex}A filtered simplicial complex is a sequence
of simplicial complexes, $\mathcal{K}=\left\{ K_{j}\right\} _{j\geq0}$,
such that\[
\widehat{C}_{n}\left(K_{0}\right)\subset\widehat{C}_{n}\left(K_{1}\right)\subset\widehat{C}_{n}\left(K_{2}\right)\subset\widehat{C}_{n}\left(K_{3}\right)\subset\cdots,\]
for all $n\geq0$, where $\widehat C_n(K)$ is the collection of all $n$-simplices in the complex $K$.
\end{definition}
We say that a simplex $\sigma$ enters the filtered complex $\mathcal{K}$
at the {\it entrance time} $i$ if $\sigma\in K_{i}$ and $\sigma\notin K_{j}$ for
all $j<i$. Occasionally, we shall use the term filtration instead of filtered
complex.

The  usual computation of the homology groups of the complexes $K_{j}$
is done for each $j$ at a time, and so does
not allow for comparison of homologies between complexes. The idea of
persistent homology is to take the filtration into account and so
be able to describe how homological properties
persist or disappear as $k$ grows.

Denoting the $n$-cycles of a complex $K$ by $Z_n(K)$ and the $n$-boundaries
by $B_n(K)$, note that any cycle in $Z_{n}(K_{j})$ also belongs to
$Z_{n}(K_{j+1})$ and boundaries in $B_{n}(K_{j})$  belong to
$B_{n}(K_{j+1})$. This allows
us to define the linear maps
\beqq
i_{*}^{\left(n,j\right)}:H_{n}\left(K_{j}\right)&\longrightarrow& H_{n}\left(K_{j+1}\right),\\
\bar{z}=z+B_{n}\left(K_{j}\right)&\longmapsto&\bar{z}=z+B_{n}\left(K_{j+1}\right).
\eeqq

\begin{definition}
The $p$-persistent $n$-th homology group of $K_{j}$ is defined
by
\beqq
H_{n}^{p}\left(K_{j}\right)=i_{*}^{p}\left(H_{n}\left(K_{j}\right)\right)\subset H_{n}\left(K_{j+p}\right),
\eeqq
where $i_{*}^{p}$  denotes the composition
 \beqq
i_{*}^{\left(n,j+p-1\right)}\circ i_{*}^{\left(n,j+p-2\right)}\circ\cdots\circ i_{*}^{\left(n,j\right)}.
\eeqq
\end{definition}
The non-zero elements of the persistent homology group are the images
of $n$-cycles which exist at time $j$ (i.e.\ belong to $K_{j}$) and which `survive'
until time $j+p$, in the sense that they are not nullified by becoming
a boundary.

We now restrict the discussion to filtered complexes of finite type.
\begin{definition}
We say that a filtered simplicial complex $\mathcal{K}=\left\{ K_{j}\right\} _{j\geq0}$
is of finite type if, for all $j\geq0$ and $n\geq0$, $\widehat{C}_{n}\left(K_{j}\right)$
is finite and if there exists an index $i$ such that $K_{j}=K_{i}$
for all $j\geq i$.
\end{definition}

Now fix $n\geq 0$. Recall that we are working with $\Z_2$, so that the
$H_{n}(K_{j})$ are all vector
spaces. Using algebraic arguments%
\footnote{
For details see \cite{Carlsson2004}.%
} it can be shown that for any filtered complex of a finite type it
is possible to choose bases $\{ c_{1}^{j},c_{2}^{j},\ldots,c_{m_{j}}^{j}\} $,
of $H_{n}(K_{j})$, for all $j\geq0$, such that for any
$1\leq k\leq m_{j}$, $i_{*}(c_{k}^{j})\in\{ 0,c_{1}^{j+1},c_{2}^{j+1},\ldots,c_{m_{j+1}}^{j+1}\} $
and $i_{*}(c_{k}^{j})=i_{*}(c_{k'}^{j})$ for
$k\neq k'$ only if $i_{*}(c_{k}^{j})=0$.

Figure 6.1 shows an example of this relation between
the bases for a certain filtered complex. The elements below each
of the homologies form a basis of the homology. Note that we have written
 $H_k\equiv H_n(K_k)$ in order to save space.
%
\beqq
\begin{array}{ccccccccccccccccc}
\label{fig:Bases-across}
  H_{0} &  & H_{1} &  & H_{2} &  & H_{3} &  & H_{4} &  & H_{5} &  & H_{6} & & \cdots\\
  c_{1}^{0} &\stackrel{i_{*}}{\mapsto}   & c_{1}^{1}& \stackrel{i_{*}}{\mapsto}   & c_{1}^{2} &\stackrel{i_{*}}{\mapsto}   &
 c_{1}^{3}&\stackrel{i_{*}}{\mapsto}   & c_{1}^{4} & \stackrel{i_{*}}{\mapsto} & 0\\
 &  & c_{2}^{1} &\stackrel{i_{*}}{\mapsto}  & c_{2}^{2} & \stackrel{i_{*}}{\mapsto}  & 0\\
  &  &  &  &  &  & c_{2}^{3} &\stackrel{i_{*}}{\mapsto} & c_{2}^{4} &\stackrel{i_{*}}{\mapsto} & c_{1}^{5} &\stackrel{i_{*}}{\mapsto} & c_{1}^{6} & \stackrel{i_{*}}{\mapsto} & \cdots\\
  &  &  &  &  &  &  & &   c_{3}^{4} & \stackrel{i_{*}}{\mapsto}  & c_{2}^{5}
&\stackrel{i_{*}}{\mapsto} & 0\\
 &  &  &  &  &  &  &  &  &  &   c_{3}^{5}& \stackrel{i_{*}}{\mapsto}   & 0\\
  &  &  &  &  &  &  &  &  &  &  &  &  c_{2}^{6}& \stackrel{i_{*}}{\mapsto} & \cdots
\end{array}
\eeqq
\bc
{{\footnotesize \smc Fig 6.1}} {\footnotesize \it Bases across the filtration.}
\ec
\vskip0.15truein

For any basis element $c_{k}^{j}$ which is not an image of a previous
basis element $c_{k}^{j-1}$, either there is a minimal number $p\geq1$
such that $i_{*}^{p}(c_{k}^{j})=0$ or $i_{*}^{p}(c_{k}^{j})\neq0$
for any $p\geq1$. Each such element can be matched
to the interval $(j,j+p)$, or, respectively, $(j,\infty)$.
These intervals give rise to the graphical presentation in the form
of a barcode. The horizontal axis of the barcode scheme represents
time - i.e.\ the index of the complex - and each bar spanning the interval
from $j$ to $i$ corresponds to an interval $(j,i)$ linked
to a basis element as above.

\setcounter{figure}{1}
\begin{figure}[h]
\begin{centering}
\includegraphics[scale=0.45]{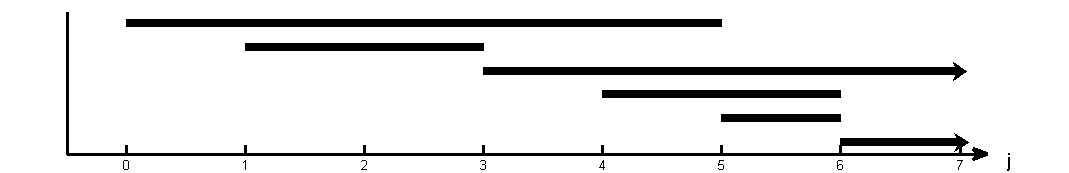}
\caption{\label{fig:Barcode-representation}Barcode representation of homology
bases.}
\end{centering}
\end{figure}

Thus, for the bases of Figure 6.1, assuming
that $K_{j}=K_{6}$ for all $j\geq6$, the barcode is as
in Figure \ref{fig:Barcode-representation}.

Note that the Betti numbers of each of the complexes in the filtration
can be easily derived from its barcode:  $\beta_{n}(K_{j})$
is the number of bars which intersect a vertical line at time
$j$, excluding those ending at time $j$.

Recall that we considered a single fixed dimension of the homology,
$n$. The collections of bars of persistent homologies of all dimensions
is called the barcode of the filtered complex (of finite type) $\mathcal{K}=\{ K_{j}\} _{j\geq0}$.

\subsection{Random filtered complexes and entrance time fields}
\label{sec:Random-Filtered-Complexes}
Random filtered complexes are the link between the `deterministic'
notion of persistent homology and the random setting.

Assuming the ubiquitous  probability space $(\Omega,\mathcal{F},P)$,
 a random filtered complex will be defined as a mapping
from $\Omega$ to some space $F$ of filtered complexes. However, allowing
$F$ to be too general makes it virtually impossible to define a meaningful
$\sigma$-algebra on it, and so we shall restrict our discussion here
to  cases in which the elements of $F$ are all sub-complexes of
 finite universal complex. Among other things, this
implies that they are of  finite type.

We begin with a trivial generalization of the definition of a filtered
simplicial complex. We now allow the set of indices of a filtered
complex to be any well-ordered infinite set in $\mathbb{R}\cup\{ -\infty\} $,
so that a filtered complex is now of the form $\mathcal{K}=\{ K_{\alpha}\} _{\alpha\in A}$.
The definitions which followed the definition of filtered complexes
can be easily adjusted accordingly. The motivation
for this is that, while in the deterministic case one deals with a
fixed filtered complex, with a fixed set of indices,
in the random case one needs to assign indices meaningful to a wider
possible  set of outcomes, and  the natural numbers no longer suffice
as an index set.

In addition, for the discussion of random filtered complexes, we think
of finite type complexes as having only a finite number of indices by
discarding the constant tail of complexes. Formally, we
can think about it as restricting ourselves to discussing only complexes
with tail defined in a canonical way: for all filtrations, $\{ K_{\alpha}\} _{\alpha\in A}$,
there exists $\alpha_{0}\in A$ such that for all $\alpha_{0}\leq\alpha\in A$,
$K_{\alpha_{0}}=K_{\alpha}$ and such that $A\cap[\alpha_{0},\infty)=\{ \alpha_{0}+i\} _{i=0}^{\infty}$.
To save some tedious notation we simply work with a finite portion
of the complex.

Next, some notation. For a given simplex $\sigma$ in a filtered complex,
we define its entrance time
\beqq
\mbox{ent}\ \equiv\ \ \mbox{ent}(\sigma) \ \definedas \
\min\left\{\alpha\: \sigma\in K_\alpha\right\},
\eeqq
if the minimum is finite and $\infty$ otherwise.
For a simplicial complex $K$, let
$\mathfrak{F}(K)$ denote all finite type filtrations,
 $\{ K_{\alpha}\} _{\alpha\in A}$, of $K$ satisfying
the condition that, for any $\alpha\in A$, there exists a simplex $\sigma\in K$
with entrance time $\mbox{ent}(\sigma)=\alpha$.
 This condition basically says that we
consider filtrations with no `spare' complexes, which, loosely
speaking, contain no additional information.

Note that we then have the natural injective mapping
\beqq
\pi \ \equiv \
\pi_{K}\:\mathfrak{F}\left(K\right)\longrightarrow\prod_{\sigma\in K}\overline{\mathbb{R}}_{\sigma}=\overline{\mathbb{R}}^{\mbox{card}\left(K\right)}\,\,\,,\,\,\,\left\{ K_{\alpha}\right\} _{\alpha\in A}\longmapsto\left\{ \mbox{ent}\left(\sigma\right)\right\} _{\sigma\in K}.
\eeqq

Using the mapping $\pi$, $\mathfrak{F}(K)$ can be endowed
with the structure of a measurable space, determined by the
rule: $F\subset\mathfrak{F}(K)$ is measurable if and only
if $F=\pi^{-1}(B)$ for some Borel set $B$ in $\prod_{\sigma\in K}\overline{\mathbb{R}}_{\sigma}$,
with respect to the standard product topology on $\prod_{\sigma\in K}\overline{\mathbb{R}}_{\sigma}$.
We denote this $\sigma$-algebra by $\mathfrak{B}(K)$.

Note that $\mathfrak{B}(K)$ is the Borel $\sigma$-algebra on $\mathfrak{F}(K)$,
when endowing it with the topology defined by the rule: $F\subset\mathfrak{F}(K)$
is open if and only if $F=\pi^{-1}(G)$ for some open set
$G$ in $\prod_{\sigma\in K}\overline{\mathbb{R}}_{\sigma}$.

We are now finally ready to define the random filtration of a complex.
\begin{definition}
Let $K$ be a fixed universal complex. A random filtration of $K$
is a measurable function $\mathcal{K}:(\Omega,\mathcal{F})\rightarrow(\mathfrak{F}(K),\mathfrak{B}(K))$.
\end{definition}

 The following lemma, in which  $\overline{\mathbb{R}}$  denotes the two point compactification of $\mathbb{R}$, is now straightforward.

\begin{Lemma}
\label{lem:filtration-measurability}For a finite complex $K$, the
mapping $\mathcal{K}:(\Omega,\mathcal{F})\rightarrow(\mathfrak{F}(K),\mathfrak{B}(K))$
is measurable if and only if $\mbox{ent}_{\mathcal{K}(\omega)}(\sigma):(\Omega,\mathcal{F})\rightarrow(\overline{\mathbb{R}},\mathfrak{B}(\overline{\mathbb{R}}))$
is measurable for all $\sigma\in K$ (where $\mathfrak{B}(\overline{\mathbb{R}})$
is the Borel $\sigma$-algebra on $\overline{\mathbb{R}}$).
\end{Lemma}
Lemma \ref{lem:filtration-measurability} implies that $\pi\circ\mathcal{K}$
is a random field on $K$. The next result shows that under certain
compatibility conditions on a field on $K$, the converse is also
true.
\begin{definition}
\label{def:Entrance-Times-Field}Let $E=\left\{ E_{\sigma}\right\} _{\sigma\in K}$
be a random field on a finite simplicial complex $K$. If $E_{\sigma}\left(\omega\right)\leq E_{\tau}\left(\omega\right)$
for any simplices $\sigma,\tau\in K$ for which $\sigma\subset\tau$,
we say that $E$ is an entrance time field on $K$.
\end{definition}
\begin{Corollary}
If $E=\left\{ E_{\sigma}\right\} _{\sigma\in K}$ is an entrance time
field on a finite simplicial complex $K$, then $\pi_{K}^{-1}\left(E\right)$
is a random filtration of $K$. Moreover, $\pi_{K}$ gives a 1-1
correspondence between random filtrations and entrance time fields
on $K$.
\end{Corollary}
Note that even when a field $E^{'}$ on $K$ does not satisfy the
condition of Definition \ref{def:Entrance-Times-Field} we can define
an entrance time field by $E=\max\{ E_{\tau}^{'}:\tau\subset\sigma\} $.


\bibliographystyle{plain}
\bibliography{hrf-bib}

\address{Robert J.\ Adler \\Electrical Engineering\\ Technion, Haifa,
Israel 32000 \\
\printead{e1}\\
\printead{u1}}

 \address{Omer Bobrowski \\Electrical Engineering\\ Technion, Haifa,
Israel 32000 \\
\printead{e2}\\
\printead{u2}}

 \address{Matthew Strom Borman \\Mathematics \\ University of Chicago
5734 S. University Ave\\
Chicago, IL 60637 \\
\printead{e4}\\
\printead{u4}}

\address{Eliran Subag\\ Electrical Engineering\\ Technion, Haifa,
Israel 32000 \\
\printead{e3}
}

\address{Shmuel Weinberger\\Mathematics \\ University of Chicago
5734 S. University Ave\\
Chicago, IL 60637 \\
\printead{e5}\\
\printead{u5}}

\end{document}